\documentclass[12pt, a4paper]{article}

\usepackage{amsmath,amsthm,amssymb}
\usepackage[truedimen,margin=25truemm]{geometry}
\usepackage{graphicx,color,psfrag}
\usepackage{mathrsfs}

\theoremstyle{definition}
\newtheorem{definition}{Definition}[section]
\newtheorem{thm}[definition]{Theorem}

\newtheorem{lem}[definition]{Lemma}
\newtheorem{cor}[definition]{Corollary}
\newtheorem{rei}[definition]{Example}
\newtheorem{remark}[definition]{Remark}

\DeclareMathOperator{\Nu}{N}

\DeclareMathOperator{\sgn}{sgn}
\DeclareMathOperator{\perm}{perm}

\title{Principal Specialization of Monomial Symmetric Polynomials and Group Determinants of Cyclic Groups}
\author{Naoya Yamaguchi, Yuka Yamaguchi and Genki Shibukawa}

\begin{document}

\maketitle

\begin{abstract}
In this paper, we study the principal specialization of monomial symmetric polynomials and investigate the special values of these polynomials at 
\[
\zeta_{(n,k)} := ( 1, \zeta_n, \zeta_n^2, \dots, \zeta_n^{kn-1} ),
\]
where \(\zeta_n\) is a primitive \(n\)th root of unity. 
We give explicit formulas for several classes of special values. 
We also show that these special values naturally appear as the coefficients in the expansion of the $k$th power of the circulant determinant of order $n$ (the group determinant of the cyclic group of order $n$). 
These results extend Ore's formulas for the case $k = 1$. 
Furthermore, we determine the number of terms in the $k$th power of the group permanent of the cyclic group of order $n$. This extends Brualdi and Newman's result for $k = 1$. 
\end{abstract}

\section{Introduction}

\subsection{Monomial symmetric polynomial (MSP)}

For a positive integer $N$, let the symmetric group $\mathfrak{S}_{N}$ of degree $N$ act on $\mathbb{Z}^{N}$ by
\begin{align}
\begin{array}{cccccc}
\mathfrak{S}_{N} & \curvearrowright & \mathbb{Z}^{N} & \xrightarrow{\simeq} & \mathbb{Z}^{N} & \nonumber \\
\rotatebox{90}{$\in$}  & & \rotatebox{90}{$\in$} & & \rotatebox{90}{$\in$} & \\
\sigma & \curvearrowright & \lambda := (\lambda_{1}, \lambda_{2}, \ldots, \lambda_{N}) & \mapsto & \sigma \cdot \lambda := (\lambda_{\sigma^{-1} (1)}, \lambda_{\sigma^{-1}(2)}, \ldots, \lambda_{\sigma^{-1} (N)}). & 
\end{array}
\end{align}
Let $\mathbb{C}[x_{1}, x_{2}, \ldots, x_{N}]$ be the ring of polynomials in $N$ independent variables $x_{1}, x_{2}, \ldots, x_{N}$ with complex coefficients, and define the action of $\mathfrak{S}_{N}$ on $\mathbb{C}[x_{1}, x_{2}, \ldots, x_{N}]$ by permutation of the variables $x_{i}\mapsto \sigma (x_{i}):=x_{\sigma (i)}$. 
We consider the subring 
$$
R_{N}
   :=
   \mathbb{C}[x_{1}, x_{2}, \ldots, x_{N}]^{\mathfrak{S}_{N}}
   =
   \{f \in \mathbb{C}[x_{1}, x_{2}, \ldots, x_{N}] \mid \sigma (f)=f \: \: \text{for any} \: \: \sigma \in \mathfrak{S}_{N} \},
$$
and we call the elements of $R_{N}$ symmetric polynomials. 
For any non-negative integer $k$, 
let 
$$
R_{N}^{k} := \{ 0 \} \cup \{ f \in R_{N} \mid f \: \: \text{is homogeneous of degree} \: k \}. 
$$
Then $R_{N}^{k}$ is a finite-dimensional vector space over $\mathbb{C}$, 
and $R_{N}$ is a graded ring with the following direct sum decomposition:
$$
R_{N}=\bigoplus_{k\geq 0}R_{N}^{k}.
$$

For $\lambda \in \mathbb{Z}^N$, 
we define
\[
m_{\lambda}(x)
:=
\sum_{\mu \in \mathfrak{S}_N \cdot \lambda} x^\mu \in \mathbb{C}[x_{1}^{\pm1}, x_{2}^{\pm1}, \ldots, x_{N}^{\pm1}]^{\mathfrak{S}_N},
\]
where
\[
\mathfrak{S}_N \cdot \lambda
:=
\{\sigma \cdot \lambda \mid \sigma \in \mathfrak{S}_N\},
\qquad
x^\mu := x_1^{\mu_1} x_2^{\mu_2} \cdots x_N^{\mu_N}.
\]
The Laurent polynomial $m_\lambda(x)$ is called the Monomial Symmetric (Laurent) Polynomial (MSP). 
We denote the set of partitions of length $\leq N$ by
\[
\mathcal{P}_{N}
:=
\left\{
(\lambda_{1}, \lambda_{2}, \ldots, \lambda_{N}) \in \mathbb{Z}^{N}
\mid
0\leq \lambda_{1} \leq \lambda_{2} \leq \cdots \leq \lambda_{N}
\right\}.
\]
If $\lambda \in \mathcal{P}_N$, 
then $m_\lambda(x)$ is the usual monomial symmetric polynomial. 

For $\lambda\in\mathbb Z^N$, 
let 
\[
\lambda[a]
:=
| \{ j \in \{1, 2, \ldots, N \}\mid \lambda_j=a\} |
\qquad (a \in \mathbb{Z}).
\]
Then the stabilizer subgroup of $\lambda$ in $\mathfrak{S}_N$ has order
\[
|\mathfrak{S}_N^\lambda|
=
\prod_{a \in \mathbb{Z}} \lambda[a]!,
\]
where only finitely many factors are different from $1$.
Hence we have
\begin{align}
m_{\lambda}(x)
=
\left(
\prod_{a \in \mathbb{Z}}\frac{1}{\lambda[a]!}
\right)
\sum_{\sigma \in \mathfrak{S}_N}x^{\sigma \cdot \lambda}.
\label{eq:MSP}
\end{align}
For each non-negative integer $k$, 
\[
\{m_\lambda(x) \mid \lambda \in \mathcal{P}_N,\ |\lambda|=k\}
\]
is the standard basis of $R_N^k$. 

For $k \geq 1$, consider the special partitions $\lambda =(k):=(0, 0, \ldots, 0, k)$. 
Then 
\[
m_{(k)}(x)
   =
   p_{k}(x)
   :=
   \sum_{i=1}^{N}x_{i}^{k}, 
\]
which is the $k$th power sum. 
Also, 
for $1 \leq k \leq N$, 
consider 
\[
\lambda =(1^{k}):=(0, 0, \ldots, 0, \overbrace{\rule{0pt}{9pt}1, 1, \ldots, 1}^{k}). 
\]
Then 
$$
m_{(1^{k})}(x)
   =
   e_{k}(x)
   :=
   \sum_{1\leq i_{1}<\cdots <i_{k}\leq N} x_{i_{1}} x_{i_{2}} \cdots x_{i_{k}}, 
$$
which is the $k$th elementary symmetric polynomial. 

One of the standard specializations of symmetric polynomials is the principal specialization
\[
x = (1, t, t^2, \ldots, t^{N-1}), \qquad t \in \mathbb{C}.
\]
In this paper, for $N=kn$, we consider the specialization
\[
x = \zeta_{(n,k)}
:=
(1, \zeta_n, \zeta_n^2, \ldots, \zeta_n^{kn-1}) \in \mathbb{C}^{kn},
\]
where $\zeta_n$ is a primitive complex $n$th root of unity.
We study the special values
\[
m_{\lambda}(\zeta_{(n, k)}).
\]
These values are related to zonal spherical functions of Gelfand pairs for complex reflection groups~\cite{MR2047846}.
Moreover, from a generating function for $m_{\lambda}(\zeta_{(n, k)})$
(see Remark~\ref{rem:4.1}), these values also appear in the generalized Waring formula~\cite{MR1401004}.
Our first main theorem gives explicit formulas for several classes of partitions $\lambda$. 

\begin{definition}\label{def:Lambda}
Let $n$ and $k$ be positive integers. 
We set
\[
\Lambda_{n}^{k} 
:= 
\left\{
(\lambda_{1}, \lambda_{2}, \ldots, \lambda_{kn}) \in \mathbb{Z}^{kn}
\mid
1 \leq \lambda_{1} \leq \lambda_{2} \leq \cdots \leq \lambda_{kn} \leq n
\right\}.
\]
For $\lambda = (\lambda_1, \lambda_{2}, \ldots, \lambda_{kn})$, 
put
\[
|\lambda| := \lambda_{1} + \lambda_{2} + \cdots + \lambda_{kn}.
\]
We also set
\[
\tilde{\Lambda}_{n}^{k}
:=
\left\{
\lambda \in \Lambda_n^k
\mid
|\lambda|\equiv 0 \! \! \pmod{n}
\right\}.
\]
\end{definition}

\begin{remark}\label{rem:1.2}
For any $a\in \mathbb{Z}$, 
let $\overline{a}$ denote the unique integer satisfying $1 \leq \overline{a} \leq n$ and $\overline{a} \equiv a \pmod{n}$.
For $\lambda = (\lambda_1, \lambda_2, \ldots, \lambda_{kn}) \in \mathbb{Z}^{kn}$, 
let $\overline{\lambda} \in \Lambda_n^k$ be the nondecreasing rearrangement of
\[
(\overline{\lambda_1}, \overline{\lambda_2}, \ldots, \overline{\lambda_{kn}}).
\]
Then
\[
| \lambda | \equiv | \overline{\lambda} |\pmod n.
\]
Moreover, the full permutation sums corresponding to $\lambda$ and
$\overline{\lambda}$ have the same value at $\zeta_{(n,k)}$. Hence
\[
m_\lambda(\zeta_{(n, k)})
=
c \, m_{\overline{\lambda}}(\zeta_{(n, k)}), 
\qquad c := \frac{
\prod_{b=1}^{n} \overline{\lambda}[b]!}{\prod_{a\in\mathbb{Z}} \lambda[a]!} \neq 0. 
\]
In particular,
\[
m_\lambda(\zeta_{(n, k)}) \neq 0 \iff 
m_{\overline{\lambda}}(\zeta_{(n,k)}) \neq 0.
\]
\end{remark}

\begin{thm}\label{thm:1.1}
Let $n$ and $k$ be positive integers. 
The following statements are true:  
\begin{enumerate}
\item[\rm (1)] For any $\lambda = (\lambda_{1}, \lambda_{2}, \overbrace{\rule{0pt}{9pt}n, n, \ldots, n}^{k n - 2}) \in \tilde{\Lambda}_{n}^{k}$ with $\lambda_{1} \leq \lambda_{2} < n$, we have 
\begin{align*}
m_{\lambda}(\zeta_{(n, k)}) = 
\begin{cases}
- \frac{n}{2} k \neq 0, & \lambda_{1} = \lambda_{2}, \\ 
- k n \neq 0, & \lambda_{1} \neq \lambda_{2}. 
\end{cases}
\end{align*}
\item[\rm (2)] For any $\lambda = (\lambda_{1}, \lambda_{1}, \lambda_{1}, \overbrace{\rule{0pt}{9pt}n, n, \ldots, n}^{k n - 3}) \in \tilde{\Lambda}_{n}^{k}$ with $\lambda_{1} < n$, 
we have 
\[
m_{\lambda}(\zeta_{(n, k)}) = \frac{n}{3} k \neq 0. 
\]
\item[\rm (3)] 
For any $\lambda \in \tilde{\Lambda}_{n}^{k}$ of either of the following forms: 
\[
\lambda =
(\lambda_{1}, \lambda_{1}, \lambda_{2},
\overbrace{\rule{0pt}{9pt}n, n, \ldots, n}^{k n - 3}),
\quad
\lambda =
(\lambda_{1}, \lambda_{2}, \lambda_{2},
\overbrace{\rule{0pt}{9pt}n, n, \ldots, n}^{k n - 3}),
\]
where $\lambda_{1} < \lambda_{2} < n$, 
we have 
\[
m_{\lambda}(\zeta_{(n, k)}) = k n \neq 0. 
\]
\item[\rm (4)] For any $\lambda = (\lambda_{1}, \lambda_{2}, \lambda_{3}, \overbrace{\rule{0pt}{9pt}n, n, \ldots, n}^{k n - 3}) \in \tilde{\Lambda}_{n}^{k}$ with $\lambda_{1} < \lambda_{2} < \lambda_{3} < n$, 
we have 
\[
m_{\lambda}(\zeta_{(n, k)}) = 2 k n \neq 0. 
\]
\item[\rm (5)] For any $\lambda = (\overbrace{\rule{0pt}{9pt}\lambda_{1}, \lambda_{1}, \ldots, \lambda_{1}}^{a}, \overbrace{\rule{0pt}{9pt}n, n, \ldots, n}^{k n - a}) \in \mathbb{Z}^{k n}$ with $n \nmid \lambda_{1}$, $n \mid a \lambda_{1}$ and $0 \leq a \leq kn$, we have 
$$
m_{\lambda}(\zeta_{(n, k)}) 
= 
(-1)^{a + \frac{a}{n} \gcd(\lambda_{1}, n)} 
\binom{k \gcd(\lambda_{1}, n)}{\frac{a}{n} \gcd(\lambda_{1}, n)} \neq 0. 
$$
\item[\rm (6)] For any $\lambda = (\overbrace{\rule{0pt}{9pt}\lambda_{1}, \lambda_{1}, \ldots, \lambda_{1}}^{a}, \overbrace{\rule{0pt}{9pt}\lambda_{2}, \lambda_{2}, \ldots, \lambda_{2}}^{k n - a}) \in \mathbb{Z}^{k n}$ with $n \nmid \lambda_{2} - \lambda_{1}$ and $1 \leq a \leq k n$, we have 
$$
m_{\lambda}(\zeta_{(n, k)}) 
= (- 1)^{k (n + 1) \lambda_{1}} m_{\lambda'}(\zeta_{(n, k)}), 
$$
where $\lambda' := (\overbrace{\rule{0pt}{9pt}\lambda_{2} - \lambda_{1}, \lambda_{2} - \lambda_{1}, \ldots, \lambda_{2} - \lambda_{1}}^{k n - a}, \overbrace{\rule{0pt}{9pt}n, n, \ldots, n}^{a})$. 
\item[\rm (7)] Let $l$ be a positive integer. 
When 
\[
\lambda = (\lambda_{1}, \lambda_{2}, \ldots, \lambda_{k n}) \in \Lambda_{n}^{k} \: \: \text{and} \: \: \mu = (\mu_{1}, \mu_{2}, \ldots, \mu_{(k + l) n}) \in \Lambda_{n}^{k + l}
\]
satisfy the condition
\begin{flalign*}
\left| \left\{ i \mid \lambda_{i} = a, 1 \leq i \leq k n \right\} \right| \leq \left| \left\{ i \mid \mu_{i} = a, 1 \leq i \leq (k + l) n \right\} \right| \:\: \text{for any}\:\: 1 \leq a \leq n, 
\end{flalign*}
we write this as $\lambda \triangleleft \mu$ and we define $\mu \backslash \lambda \in \Lambda_{n}^{l}$ as the sequence obtained by removing the $\lambda_{i}$ from $\mu$. 
Then, for any $\mu \in \Lambda_{n}^{k + l}$, we have 
$$
m_{\mu}(\zeta_{(n, k+l)}) = \displaystyle\sum_{\substack{\lambda \in \Lambda_{n}^{k} \\ \lambda \: \triangleleft \: \mu}} m_{\lambda}(\zeta_{(n, k)}) m_{\mu \backslash \lambda}(\zeta_{(n, l)}). 
$$
\end{enumerate}
\end{thm}

The cases where $k = 1$ and $\lambda_{i}$ are different from each other in Theorem~$\ref{thm:1.1}$~(1) and (4) are given in \cite[p.~345]{MR42365}, 
and the case of $k = 1$ in Theorem~$\ref{thm:1.1}$~(5) is given in \cite[p.~346, Eq~(18)]{MR42365}. 
All assertions except (5) and (6) are proved by using our second main theorem (Theorem~$\ref{thm:1.2}$) on the circulant determinant.

From the point of view of the circulant determinant, the above results on $m_{\lambda}(\zeta_{(n, k)})$ describe the non-vanishing properties of some terms in the $k$th power of the circulant determinant of order $n$. 

\begin{remark}
Riordan~\cite{MR281628} and Merca~\cite{Merca2015, MR3742496} studied expansions of augmented monomial symmetric polynomials in terms of power sums $m_{(k)}(x) = p_{k}(x)$. 
Their results give other expressions for the special values considered in this paper. 
In the present paper, 
however, 
we keep the circulant determinant approach, 
because our main purpose is to relate these special values to powers of circulant determinants and group determinants.
\end{remark}

\subsection{Number of terms in the group permanent and the group determinant of a cyclic group}\label{sec1.2}

For a positive integer $n$, 
the determinant 
$$
{\rm C}(x_{1}, x_{2}, \ldots, x_{n}) := \det{\begin{pmatrix} x_{n} & x_{n - 1} &  \cdots & x_{1} \\ x_{1} & x_{n} & \cdots & x_{2} \\ \vdots & \vdots & \ddots & \vdots \\ x_{n - 1} & x_{n - 2} & \cdots & x_{n} \end{pmatrix}}
$$
is called the circulant determinant of order $n$. 
It is well known that the circulant determinant can be factorized over $\mathbb{C}$ as follows: 
\begin{align}\label{eq:circulant-factorization}
{\rm C}(x_{1}, x_{2}, \ldots, x_{n}) = \prod_{i = 1}^{n} \sum_{j = 1}^{n} \zeta_{n}^{ij} x_{j}. 
\end{align}
By expanding this factorization, 
Ore~\cite[THEOREM~1]{MR42365} obtained the explicit formula for the circulant determinant: 
\begin{align*}
{\rm C}(x_{1}, x_{2}, \ldots, x_{n})  = \sum_{\lambda \in \Lambda_{n}^{1}} m_{\lambda}(\zeta_{(n, 1)}) x_{\lambda} = \sum_{\lambda \in \tilde{\Lambda}_{n}^{1}} m_{\lambda}(\zeta_{(n, 1)}) x_{\lambda}, 
\end{align*}
where $x_{\lambda} := x_{\lambda_1} x_{\lambda_2} \cdots x_{\lambda_n}$.
We present the explicit formula for a power of the circulant determinant. 

\begin{thm}\label{thm:1.2}
For any positive integer $k$, 
we have 
\begin{align*}
{\rm C}(x_{1}, x_{2}, \ldots, x_{n})^{k} = \sum_{\lambda \in \Lambda_{n}^{k}} m_{\lambda}(\zeta_{(n, k)}) x_{\lambda} = \sum_{\lambda \in \tilde{\Lambda}_{n}^{k}} m_{\lambda}(\zeta_{(n, k)}) x_{\lambda}, 
\end{align*}
where $x_{\lambda} := x_{\lambda_1} x_{\lambda_2} \cdots x_{\lambda_{kn}}$. 
\end{thm}

The group determinant is a generalization of the circulant determinant. 
For a finite group $G := \{ g_{1}, g_{2}, \ldots, g_{n} \}$, 
let $x_{g}$ be an indeterminate for each $g \in G$, 
and let $\mathbb{Z}[x_{g}] := \mathbb{Z}[x_{g_{1}}, x_{g_{2}}, \ldots, x_{g_{n}}]$ be the multivariate polynomial ring in the $x_{g}$ over $\mathbb{Z}$. 
The group determinant $\Theta(G)$ of $G$ is defined as follows (see e.g., \cite{MR1085397, johnson2019group}): 
$$
\Theta(G) := \det{(x_{g h^{-1}})_{g, h \in G}} = \sum_{\sigma \in \mathfrak{S}_{n}} \sgn(\sigma) x_{g_{1} g_{\sigma(1)}^{-1 }} x_{g_{2} g_{\sigma(2)}^{- 1}} \cdots x_{g_{n} g_{\sigma(n)}^{- 1}} \in \mathbb{Z}[x_{g}]. 
$$
From this definition, 
it follows that $\Theta(G)$ is a homogeneous polynomial of degree $n$ in $n$ variables. 
When $G$ is abelian, 
for any term $x_{a_{1}} x_{a_{2}} \cdots x_{a_{n}}$ in $\Theta(G)$, 
the product $a_{1} a_{2} \cdots a_{n}$ is the unit element $e$ of $G$. 
An analogous statement also holds for non-abelian groups after choosing the appropriate order of the factors; see \cite[Lemma~1]{MR1123661}. 
Note that if $G$ is the cyclic group $\mathbb{Z} / n \mathbb{Z}$, 
then $\Theta(G)$ equals the circulant determinant of order $n$. 

For a finite group $G$, the group permanent ${\rm P}(G)$ of $G$ is defined by 
\[
{\rm P}(G) := \perm(x_{g h^{-1}})_{g, h \in G} = \sum_{\sigma \in \mathfrak{S}_n} x_{g_{1} g_{\sigma(1)}^{-1 }} x_{g_{2} g_{\sigma(2)}^{- 1}} \cdots x_{g_{n} g_{\sigma(n)}^{- 1}} \in \mathbb{Z}[x_{g}]. 
\]
For any $f(x_{g}) \in \mathbb{Z}[x_{g}]$, 
let $\Nu(f(x_{g}))$ denote the number of terms in $f(x_{g})$. 
Hall~\cite{MR50579} proved that when $G$ is an abelian group, 
for any non-negative integers $i_{1}, i_{2}, \ldots, i_{n}$ with $i_{1} + i_{2} + \cdots + i_{n} = n$ and $g_{1}^{i_{1}} g_{2}^{i_{2}} \cdots g_{n}^{i_{n}} = e$, 
there exists a permutation $\sigma \in \mathfrak{S}_{n}$ satisfying $x_{g_{1} g_{\sigma(1)}^{-1}} x_{g_{2} g_{\sigma(2)}^{-1}} \cdots x_{g_{n} g_{\sigma(n)}^{-1}} = x_{g_{1}}^{i_{1}} x_{g_{2}}^{i_{2}} \cdots x_{g_{n}}^{i_{n}}$. 
From this, we immediately obtain the following \cite[p.~121]{MR2747803}: If $G$ is an abelian group, then
\[
\Nu({\rm P}(G)) = \left| \left\{ (i_{1}, i_{2}, \ldots, i_{n}) \in \mathbb{Z}_{\geq 0}^{n} \mid i_{1} + i_{2} + \cdots + i_{n} = n, \: g_{1}^{i_{1}} g_{2}^{i_{2}} \cdots g_{n}^{i_{n}} = e \:  \right\} \right|. 
\]
Let $\varphi$ be Euler's totient function, 
and let 
\[
S(n, m) := \{ (i_{1}, i_{2}, \ldots, i_{n}) \in \mathbb{Z}_{\geq 0}^{n} \mid i_{1} + i_{2} + \cdots + i_{n} = m, \: \: i_{1} + 2 i_{2} + \cdots + n i_{n} \equiv  0 \!\! \pmod{n} \}. 
\]
Brualdi and Newman~\cite{MR266850} showed that the following holds: 
\[
\Nu({\rm P}(\mathbb{Z} / n \mathbb{Z})) = |S(n, n)| = \frac{1}{n} \sum_{d \mid n} \binom{2d - 1}{d} \varphi\left( \frac{n}{d} \right). 
\]
We provide an explicit formula for the number of terms in the $k$th power of ${\rm P}(\mathbb{Z} / n \mathbb{Z})$. 

\begin{thm}\label{thm:1.3}
For any positive integer $n$ and $k$, we have 
\begin{align*}
\Nu({\rm P}(\mathbb{Z} / n \mathbb{Z})^{k}) 
= |\tilde{\Lambda}_{n}^{k}| 
= |S(n, kn)| 
= \frac{1}{n} \sum_{d \mid n} \binom{dk + d - 1}{d - 1} \varphi\left( \frac{n}{d} \right). 
\end{align*}
\end{thm}

Since the support of $\Theta(\mathbb{Z}/n\mathbb{Z})$ is contained in the support of
${\rm P}(\mathbb{Z}/n\mathbb{Z})$, it follows immediately that
\[
\Nu(\Theta(\mathbb{Z}/n\mathbb{Z})^{k})
\leq
\Nu({\rm P}(\mathbb{Z}/n\mathbb{Z})^{k}).
\]

By using E\u{g}ecio\u{g}lu and Remmel's~\cite{MR1137989} result on the MSP, Thomas~\cite{MR2104820} showed that, 
when $n$ is a prime power, 
$\Nu(\Theta(\mathbb{Z} / n \mathbb{Z})) = \Nu({\rm P}(\mathbb{Z} / n \mathbb{Z}))$ holds. 
Also, Colarte, Mezzetti, Mir\'o-Roig and Salat~\cite{MR3894894}, 
by using Malenfant's~\cite{malenfant2015matrixelementexpansioncirculantdeterminant} result on the circulant determinant, 
proved that if 
$\Nu(\Theta(\mathbb{Z} / n \mathbb{Z})) = \Nu({\rm P}(\mathbb{Z} / n \mathbb{Z}))$, 
then $n$ must be a prime power.  
A prime number $p$ satisfying the congruence  
\[
\binom{2p - 1}{p - 1} \equiv 1 \pmod{p^{4}}
\]
is called a Wolstenholme prime~\cite[p.~385]{MR1339137}. 
Recently, the following result was obtained~\cite[Proposition~2.3]{MR4816571}: 
For a prime $p \geq 5$ and integers $k, l \geq 1$, the congruence $| \tilde{\Lambda}_{p^{l}}^{k} | \equiv 1 \pmod{p^{2}}$ holds.
In addition, 
if $p$ is a Wolstenholme prime, 
then we have $| \tilde{\Lambda}_{p^{l}}^{1} | \equiv 1 \pmod{p^{3}}$.
Therefore, from Theorem~$\ref{thm:1.3}$, we have the following corollary. 

\begin{cor}
For prime $p \geq 5$ and integers $k, l \geq 1$, 
the congruence 
\[
\Nu({\rm P}(\mathbb{Z} / p^{l} \mathbb{Z})^{k}) 
= |\tilde{\Lambda}_{p^{l}}^{k}| \equiv 1 \pmod{p^{2}} 
\]
holds. 
In addition, 
if $p$ is a Wolstenholme prime, 
then we have 
\[
\Nu(\Theta(\mathbb{Z} / p^{l} \mathbb{Z})) 
= \Nu({\rm P}(\mathbb{Z} / p^{l} \mathbb{Z})) 
= |\tilde{\Lambda}_{p^{l}}^{1}| 
\equiv 1 \pmod{p^{3}}. 
\]
\end{cor}

Recent studies~\cite{MR4261108, MR3448603} provided an explicit formula for $\Nu({\rm P}(G))$ for any finite abelian group. 
Using the explicit formula for $\Nu({\rm P}(G))$, it is proved in \cite{MR4826770} that, for any finite abelian groups $G$ and $H$, the equality $\Nu({\rm P}(G)) = \Nu({\rm P}(H))$ holds if and only if $G \cong H$.  
This result raises the following open problem: is it true that $\Nu(\Theta(G)) = \Nu(\Theta(H))$ holds if and only if $G \cong H$? 
In \cite{MR4816571, MR4593070}, it is conjectured that, for any integer $k \geq 1$, the necessary and sufficient condition for $\Nu(\Theta(\mathbb{Z} / n \mathbb{Z})^{k}) = \Nu({\rm P}(\mathbb{Z} / n \mathbb{Z})^{k})$ to hold is that $n$ is a prime power. 
This conjecture also remains an open problem.

\section{Proofs of Theorems~$\mathbf{\ref{thm:1.2}}$ and $\mathbf{\ref{thm:1.3}}$}\label{sec2}
We prove Theorems~$\ref{thm:1.2}$ and $\ref{thm:1.3}$. 
Note that there exist integers $c_{\lambda}$ satisfying 
\begin{align*}
\Theta(\mathbb{Z} / n \mathbb{Z}) = \sum_{\lambda \in \tilde{\Lambda}_{n}^{1}} c_{\lambda} x_{\lambda}
\end{align*}
since, for a cyclic group $\mathbb{Z} / n \mathbb{Z} = \{ 1, 2, \ldots, n \}$, 
each term in $\Theta(\mathbb{Z} / n \mathbb{Z})$ is of the form $x_{i_{1}} x_{i_{2}} \cdots x_{i_{n}}$ with $n \mid i_{1} + i_{2} + \cdots + i_{n}$.

\begin{rei}\label{rei:2.1}
The group determinant of $\mathbb{Z} / 3 \mathbb{Z}$ is $x_{1}^{3} + x_{2}^{3} + x_{3}^{3} - 3 x_{1} x_{2} x_{3}$. 
The terms $x_{1}^{3}$, $x_{2}^{3}$, $x_{3}^{3}$, and $x_{1} x_{2} x_{3}$ correspond to the partitions $(1, 1, 1)$, $(2, 2, 2)$, $(3, 3, 3)$, and $(1, 2, 3)$ in $\tilde{\Lambda}_{3}^{1}$, respectively. 
\end{rei}

More generally, for the $k$th power of $\Theta(\mathbb{Z} / n \mathbb{Z})$, 
there exist integers $c_{\lambda}$ satisfying 
\begin{align*}
\Theta \left( \mathbb{Z} / n \mathbb{Z} \right)^{k} = \sum_{\lambda \in \tilde{\Lambda}_{n}^{k}} c_{\lambda} x_{\lambda}. 
\end{align*}
Theorem~$\ref{thm:1.2}$ implies that the coefficient $c_{\lambda}$ is equal to $m_{\lambda}(\zeta_{(n, k)})$. 

\begin{proof}[Proof of Theorem~$\ref{thm:1.2}$]
From the factorization~\eqref{eq:circulant-factorization} of the circulant determinant and the expression~\eqref{eq:MSP} of the MSP, we have  
\begin{align*}
\Theta \left(\mathbb{Z} / n \mathbb{Z} \right)^{k} 
&=\left(\prod_{i = 1}^{n} \sum_{j = 1}^{n} \zeta_{n}^{i j} x_{j} \right)^{k} \\
&= \left(\prod_{i = 1}^{n} \sum_{j = 1}^{n} \zeta_{n}^{i j} x_{j} \right) 
\left(\prod_{i = n + 1}^{2n} \sum_{j = 1}^{n} \zeta_{n}^{i j} x_{j} \right) \cdots 
\left(\prod_{i = (k - 1)n + 1}^{kn} \sum_{j = 1}^{n} \zeta_{n}^{i j} x_{j} \right) \\
&= \sum_{\lambda \in \Lambda_{n}^{k}} 
\left\{\sum_{\sigma \in \mathfrak{S}_{kn}} \zeta_{(n, k)}^{\sigma \cdot \lambda} \prod_{i \in \mathbb{N}} \frac{1}{\lambda[i] !} \right\} x_{\lambda} \\
&= \sum_{\lambda \in \Lambda_{n}^{k}} m_{\lambda}(\zeta_{(n, k)}) x_{\lambda}. 
\end{align*}
Moreover, by the definition of the group determinant, 
the monomials $x_{\lambda}$ with $|\lambda| \not\equiv 0 \pmod{n}$ do not appear in $\Theta(\mathbb{Z} / n \mathbb{Z})^{k}$. 
That is, $m_{\lambda}(\zeta_{(n, k)}) = 0$ holds for any $\lambda$ with $|\lambda| \not\equiv 0 \pmod{n}$.
\end{proof}


To prove Theorem~$\ref{thm:1.3}$, 
we use the following two lemmas. 

\begin{lem}[\cite{MR1691428, MR357132}]\label{lem:2.2}
It holds that
$$
|S(n, m)| = \frac{1}{n + m} \sum_{d \mid \gcd(n, m)} \binom{(n + m) / d}{n / d} \varphi(d). 
$$
\end{lem}

\begin{lem}[\cite{MR3618568}]\label{lem:2.3}
Each sequence of $2n-1$ integers contains a subsequence of $n$ terms
whose sum is divisible by $n$. 
\end{lem}

\begin{proof}[Proof of Theorem~$\ref{thm:1.3}$]
From Lemma~$\ref{lem:2.2}$, we have 
\begin{align*}
|S(n, kn)| 
&= \frac{1}{n (k + 1)} \sum_{d \mid n} \binom{(k + 1) n / d}{n / d} \varphi(d) \\ 
&= \frac{1}{n (k + 1)} \sum_{d \mid n} \binom{dk + d}{d} \varphi \left( \frac{n}{d} \right) \\
&=\frac{1}{n} \sum_{d \mid n} \binom{d k + d - 1}{d - 1} \varphi \left( \frac{n}{d} \right). 
\end{align*}
Also, from the definitions, we obtain 
\[
\Nu({\rm P}(\mathbb{Z} / n \mathbb{Z})^{k}) \leq |\tilde{\Lambda}_{n}^{k}| = |S(n, kn)|.
\]
Thus, it remains to show that, for any 
\(\lambda \in \tilde{\Lambda}_{n}^{k}\), 
the monomial \(x_{\lambda}\) appears in \({\rm P}(\mathbb{Z} / n \mathbb{Z})^{k}\). 
Let \(\lambda = (\lambda_{1}, \lambda_{2}, \ldots, \lambda_{kn}) \in \tilde{\Lambda}_{n}^{k}\). 
By using Lemma~$\ref{lem:2.3}$ iteratively, 
we find that there exists a permutation \(\sigma \in \mathfrak{S}_{kn}\) satisfying 
\begin{align*}
\lambda_{\sigma(1)} + \lambda_{\sigma(2)} + \cdots + \lambda_{\sigma(n)} 
&\equiv \lambda_{\sigma(n+1)} + \lambda_{\sigma(n+2)} + \cdots + \lambda_{\sigma(2n)} \\
&\:\;\vdots \\
&\equiv \lambda_{\sigma((k-2)n+1)} + \lambda_{\sigma((k-2)n+2)} + \cdots + \lambda_{\sigma((k-1)n)} \\
&\equiv 0 \pmod{n}.
\end{align*}
This $\sigma$ also satisfies $\lambda_{\sigma((k-1)n+1)} + \lambda_{\sigma((k-1)n+2)} + \cdots + \lambda_{\sigma(kn)} \equiv 0 \pmod{n}$ since $\lambda \in \tilde{\Lambda}_{n}^{k}$. 
For any $0 \leq i \leq k - 1$, let $\lambda^{(i)} \in \tilde{\Lambda}_{n}^{1}$ be the partition obtained by arranging the integers $\lambda_{\sigma(i n + 1)}, \lambda_{\sigma(i n + 2)}, \ldots, \lambda_{\sigma((i + 1)n)}$ in ascending order. 
From Brualdi and Newman's~\cite{MR266850} result 
$\Nu({\rm P}(\mathbb{Z} / n \mathbb{Z})) = |S(n, n)| = |\tilde{\Lambda}_{n}^{1}|$, 
it follows that the monomial $x_{\lambda^{(i)}}$ appears in ${\rm P}(\mathbb{Z} / n \mathbb{Z})$ for any $0 \leq i \leq k - 1$. 
Therefore, 
the monomial $x_{\lambda} = \prod_{i = 0}^{k - 1} x_{\lambda^{(i)}}$ appears in \({\rm P}(\mathbb{Z} / n \mathbb{Z})^{k}\) since the coefficient of every term in ${\rm P}(\mathbb{Z} / n \mathbb{Z})$ is a positive integer. 
\end{proof}

\section{Proof of Theorem~$\mathbf{\ref{thm:1.1}}$}\label{sec3}
We prove Theorem~$\ref{thm:1.1}$. 
To prove Theorem~$\ref{thm:1.1}$,  we use the following lemma. 

\begin{lem}[{\cite[Proofs of Lemmas~2 and 3]{MR1123661}}, {\cite[Lemma~3.3]{MR4227663}}]\label{lem:3.1}
Let $G$ be a finite group, 
let $e$ be the unit element of $G$ and let $n$ be the order of $G$. 
\begin{enumerate}
\item[$(1)$] If neither $a$ nor $b$ is $e$ and the monomial $x_{e}^{n-2} x_{a} x_{b}$ occurs in $\Theta(G)$, 
the coefficient of the monomial is $- n / 2$ or $- n$ depending on whether or not $a = b$. 
\item[$(2)$] If none of $a$, $b$, and $c$ is $e$ and the monomial $x_{e}^{n-3} x_{a} x_{b} x_{c}$ occurs in $\Theta(G)$, 
the coefficient of the monomial is 
\begin{enumerate}
\item[(i)] $n / 3$ if $a = b = c;$ 
\item[(ii)] $n$ if two of $a$, $b$, $c$ are equal$;$ 
\item[(iii)] $n$ if no two of them are equal and $a b \neq b a;$ 
\item[(iv)] $2n$ if no two of them are equal and $a b = b a$. $($Note that if $a b c = e$, then $a b = b a$ if and only if $a$, $b$, and $c$ commute with one another). 
\end{enumerate}
\end{enumerate}
Here, we say that a monomial occurs in a polynomial if the monomial is not canceled after combining like terms. 
\end{lem}

\begin{proof}[Proof of Theorem~$\ref{thm:1.1}$]
We first prove assertions~(1)--(4). 
To simplify the notation, 
let ${\rm C} := {\rm C}(x_{1}, x_{2}, \ldots, x_{n})$. 
Note that, from Ore's~\cite{MR42365} result ${\rm C} = \sum_{\lambda \in \tilde{\Lambda}_{n}^{1}} m_{\lambda}(\zeta_{(n, 1)}) x_{\lambda}$, the monomials $x_{\lambda}$ with $|\lambda| \not\equiv 0 \pmod{n}$ do not appear in ${\rm C}$. 
Also, 
it follows from the definition of ${\rm C}$ that the coefficient of the term $x_{n}^{n}$ in ${\rm C}$ is $1$. 
We prove (1). 
Let $\lambda = (\lambda_{1}, \lambda_{2}, n, n, \ldots, n) \in \tilde{\Lambda}_{n}^{k}$ with $\lambda_{1} \leq \lambda_{2} < n$. 
Then, the monomial $x_{n}^{n - 1} x_{\lambda_1}$ does not appear in ${\rm C}$. 
Hence the term $x_{n}^{kn - 2} x_{\lambda_1} x_{\lambda_2}$ in ${\rm C}^{k}$ can arise only by taking the term $x_{n}^{n - 2} x_{\lambda_1} x_{\lambda_2}$ from one factor of ${\rm C}$ and the term $x_n^n$ from each of the remaining $k-1$ factors. 
From Lemma~$\ref{lem:3.1}$~(1), 
we find that the coefficient of $x_{n}^{n - 2} x_{\lambda_1} x_{\lambda_2}$ in ${\rm C}$ is $-n / 2$ if $\lambda_{1} = \lambda_{2}$; $-n$ if $\lambda_{1} \neq \lambda_{2}$. 
Therefore, from Theorem~$\ref{thm:1.2}$, 
it holds that $m_{\lambda}(\zeta_{(n, k)})$ equals $-kn / 2$ if $\lambda_{1} = \lambda_{2}$; $-kn$ if $\lambda_{1} \neq \lambda_{2}$. 
In the same way, 
we can prove (2)--(4) by using Theorem~$\ref{thm:1.2}$ and Lemma~$\ref{lem:3.1}$~(2). 
Next, we prove (5). 
For this purpose, 
we consider a generating function for $m_{\lambda}(\zeta_{(n, k)})$ with $\lambda = (\overbrace{\rule{0pt}{9pt}\lambda_{1}, \lambda_{1}, \ldots, \lambda_{1}}^{a}, \overbrace{\rule{0pt}{9pt}n, n, \ldots, n}^{k n - a}) \in \mathbb{Z}^{k n}$. 
Let $u$ be an indeterminate. 
Then, we have 
\begin{align*}
\sum_{a= 0}^{kn} (-1)^{a} m_{\lambda}(\zeta_{(n, k)}) u^{k n - a} 
= \sum_{a = 0}^{kn} \sum_{1 \leq i_{1} < i_{2} < \cdots < i_{a} \leq kn} (-1)^{a} \zeta_{n}^{\lambda_{1}(i_{1} + i_{2} + \cdots + i_{a})} u^{k n - a} 
= \prod_{i = 1}^{k n } \left( u - \zeta_{n}^{i \lambda_{1}} \right). 
\end{align*}
Since $\zeta_{n}^{\lambda_{1}}$ is a primitive $l$th root of unity, 
where $d := \gcd(\lambda_{1}, n)$ and $l := \frac{n}{d}$, 
we have 
\begin{align*}
\left( u - \zeta_{n}^{\lambda_{1}} \right) \left( u - \zeta_{n}^{2 \lambda_{1}} \right) \cdots \left( u - \zeta_{n}^{l \lambda_{1}} \right) 
= \left( u - \zeta_{l} \right) \left( u - \zeta_{l}^{2} \right) \cdots \left( u - \zeta_{l}^{l} \right) 
= u^{l} - 1. 
\end{align*} 
Thus, it holds that 
\begin{align*}
\prod_{i = 1}^{k n} \left( u - \zeta_{n}^{i \lambda_{1}} \right) 
&= \prod_{i = 1}^{k d l} \left( u - \zeta_{n}^{i \lambda_{1}} \right) \\ 
&= \prod_{i = 1}^{l} \left( u - \zeta_{n}^{i \lambda_{1}} \right)^{k d} \\ 
&= (u^{l} - 1)^{kd} \\
&= \sum_{i = 0}^{kd} \binom{kd}{i} (-1)^{i} u^{l (kd - i)} \\
&= \sum_{i = 0}^{kd} \binom{kd}{i} (-1)^{i} u^{kn - \frac{n}{d} i}. 
\end{align*}
By comparing the coefficients of $u^{kn - a}$, where $\frac{n}{d} \mid a$, in the first and the last expressions, we obtain 
\begin{align*}
(-1)^{a} m_{\lambda}(\zeta_{(n, k)}) = 
\displaystyle \binom{kd}{\frac{ad}{n}} (-1)^{\frac{ad}{n}}. 
\end{align*}
Since $\frac{n}{d} \mid a$ holds when $n \mid a \lambda_{1}$, the proof of (5) is complete. 
We prove (6) by direct calculation. 
Let $[kn] := \{ 1, 2, \ldots, kn \}$ and $I^{c}:=[kn] \setminus I$. 
Then, for any 
$$
\lambda = (\overbrace{\rule{0pt}{9pt}\lambda_{1}, \lambda_{1}, \ldots, \lambda_{1}}^{a}, \overbrace{\rule{0pt}{9pt}\lambda_{2}, \lambda_{2}, \ldots, \lambda_{2}}^{k n - a}) \in \mathbb{Z}^{k n}, 
$$ 
we have
\begin{align*}
m_{\lambda}(\zeta_{(n, k)}) 
&= \sum_{\substack{I\subset [kn] \\ |I|=a}} 
\left(\prod_{i \in I} \zeta_{n}^{\lambda_{1} i} \right) \left(\prod_{j \in I^{c}} \zeta_{n}^{\lambda_{2} j} \right), \\
&= \sum_{\substack{I\subset [kn] \\ |I|=a}} 
\left(\prod_{i \in I} \zeta_{n}^{\lambda_{1} i} \right) \left(\prod_{j \in I^{c}} \zeta_{n}^{\lambda_{1} j} \right) 
\left(\prod_{j \in I^{c}} \zeta_{n}^{- \lambda_{1} j} \right)  \left(\prod_{j \in I^{c}} \zeta_{n}^{\lambda_{2} j} \right) \\
&= \sum_{\substack{I \subset [kn] \\ |I|=a}} 
\left(\prod_{i \in [kn]} \zeta_{n}^{\lambda_{1} i} \right) \left(\prod_{j \in I^{c}} \zeta_{n}^{(\lambda_{2} - \lambda_{1}) j} \right) \\
&= \zeta_{n}^{\lambda_{1} \frac{k n (k n + 1)}{2}}\sum_{\substack{I\subset [kn] \\ |I|=a}} \left(\prod_{j \in I^{c}} \zeta_{n}^{(\lambda_{2} - \lambda_{1}) j} \right). 
\end{align*}
Thus, from 
\begin{align*}
\zeta_{n}^{\frac{k n (k n + 1)}{2}} 
= 
\begin{cases}
(- 1)^{k}, & n \ \text{is even}, \\ 
1, & n \ \text{is odd}
\end{cases} 
= (- 1)^{k (n + 1)}, 
\end{align*}
we obtain 
$$
m_{\lambda}(\zeta_{(n, k)}) 
= (- 1)^{k (n + 1) \lambda_{1}} m_{\lambda'}(\zeta_{(n, k)}), 
$$ 
where $\lambda' := (\overbrace{\rule{0pt}{9pt}\lambda_{2} - \lambda_{1}, \lambda_{2} - \lambda_{1}, \ldots, \lambda_{2} - \lambda_{1}}^{k n - a}, \overbrace{\rule{0pt}{9pt}n, n, \ldots, n}^{a})$. 
Finally, from Theorem~$\ref{thm:1.2}$ and the equality $\Theta(\mathbb{Z} / n \mathbb{Z})^{k + l} = \Theta(\mathbb{Z} / n \mathbb{Z})^{k} \Theta(\mathbb{Z} / n \mathbb{Z})^{l}$, 
we obtain $(7)$. 
\end{proof}

\section{Remarks on the special values of the MSP}\label{sec4}
We give some remarks on $m_{\lambda}(\zeta_{(n, k)})$. 

\begin{remark}\label{rem:4.1}
We mention the dual Cauchy kernel formula {\rm \cite[Chapter~I~($4.2^{\prime}$)]{MR1354144}}. 
For any partitions $\lambda = (\lambda_{1}, \lambda_{2}, \ldots, \lambda_{N})$, $\mu = (\mu_{1}, \mu_{2}, \ldots, \mu_{N}) \in \mathcal{P}_{N}$, 
let $\lambda \subseteq \mu$ be the inclusion partial order defined by
$$
\lambda \subseteq \mu \quad \Longleftrightarrow \quad \lambda _{i}\leq \mu _{i}, \quad i=1, 2, \ldots, N.
$$
For positive integers $M$ and $N$, the following identity holds:
\begin{align*}
\prod_{i=1}^{M}\prod_{j=1}^{N}(1+x_{i}y_{j})
= \sum_{\lambda \subseteq (M^{N})} e_{\lambda }(x)m_{\lambda }(y),
\end{align*}
where $(M^{N}) := (\overbrace{\rule{0pt}{9pt} M, M, \ldots , M}^{N})$ and 
$e_{\lambda }(x):=\prod_{i=1}^{N}e_{\lambda _{i}}(x)$. 
As a corollary of this famous result, 
we immediately obtain the following generating function for $m_{\lambda }(\zeta_{(n, k)})$: 
For any positive integers $k$ and $n$, we have
\begin{align*}
\prod_{i=1}^{n}\prod_{j=1}^{kn}(1-x_{i} \zeta_{n}^{j-1})
   =
   \prod_{i=1}^{n}(1-x_{i}^{n})^{k} \nonumber 
   =
   \sum_{\substack{\lambda \subseteq (n^{k n}) \\ |\lambda| \equiv 0 \!\!\!\pmod{n}}}
   (-1)^{|\lambda |}e_{\lambda }(x)m_{\lambda }(\zeta_{(n, k)}).
\end{align*}
\end{remark}

\begin{remark}\label{rem:4.2}
In connection with Theorem~$\ref{thm:1.1}$~(7), 
we recall that a similar identity follows from the branching formula for Hall-Littlewood functions {\rm \cite[{Chapter}~III ($5.5{}^{\prime}$)]{MR1354144}:} 
$$
m_{\mu}(\zeta_{(n, k+l)}) = \displaystyle\sum_{\substack{\lambda \subseteq \mu \\ |\lambda| \equiv 0 \!\!\!\pmod{n}}} m_{\lambda}(\zeta_{(n, k)}) m_{\mu \backslash \lambda}(\zeta_{(n, l)}),
$$
where \(m_{\mu \backslash \lambda}\) denotes the skew monomial symmetric function
associated with the skew Young diagram \(\mu \backslash \lambda\). 
Theorem~\ref{thm:1.1}~(7), however, gives a more explicit decomposition adapted to the present specialization, 
and we do not derive it from the above branching formula in this paper. 
\end{remark}

\begin{remark}\label{rem:4.3}
The following statements are true: 
\begin{enumerate}
\item[\rm (1)] For any $\lambda \in \mathbb{Z}^{k n}$, we have $m_{\lambda}(\zeta_{(n, k)}) \in \mathbb{Z}$. 
\item[\rm (2)] For any $\lambda \in \mathbb{Z}^{k n}$ with $|\lambda| \not\equiv 0 \pmod{n}$, we have $m_{\lambda}(\zeta_{(n, k)}) = 0$.
\item[\rm (3)] For any $\lambda \in \mathbb{Z}^{k n}$ and $l \in \mathbb{Z} \setminus \{ 0 \}$ with $\gcd(l, n) = 1$, 
we have $m_{l \lambda}(\zeta_{(n, k)}) = m_{\lambda}(\zeta_{(n, k)})$, 
where $l \lambda := (l \lambda_{1}, l \lambda_{2}, \ldots, l \lambda_{k n})$. 
\end{enumerate}
These properties follow from general facts about symmetric Laurent polynomials
evaluated at roots of unity. 
In fact, (1) and (3) are true for any symmetric Laurent polynomial $f(x_{1}, x_{2}, \ldots, x_{kn})$ with integer coefficients, not just for $m_{\lambda}(\zeta_{(n, k)})$. 
That is, $f(\zeta_{(n, k)}) \in \mathbb{Z}$ and $f(1, \zeta_{n}^{l}, \ldots, (\zeta_{n}^{kn-1})^{l}) = f(\zeta_{(n, k)})$ are valid. 
Also, (2) is true for any symmetric Laurent polynomial $f(x_{1}, x_{2}, \ldots, x_{kn})$ that is homogeneous of degree $h$ with $n \nmid h$ since $f(\zeta_{(n, k)}) = f(\zeta_{n}, \zeta_{n}^{2}, \ldots, \zeta_{n}^{kn-1}, 1) = f(\zeta_{n} \cdot 1, \zeta_{n} \cdot \zeta_{n}, \ldots, \zeta_{n} \cdot \zeta_{n}^{kn-1}) = \zeta_{n}^{h} f(\zeta_{(n, k)})$ holds. 
We can also obtain (1)--(3) from Theorem~$\ref{thm:1.2}$. 
As immediate consequences of the theorem, we have (1) and (2). 
Also, (3) follows from the theorem and the fact that 
$$
\Theta(G) = \det{\left( x_{g h^{- 1}} \right)_{g, h \in G}} = \det{\left( x_{\psi(g) \psi\left(h^{- 1} \right)} \right)_{g, h \in G}} = \det{\left( x_{\psi\left( g h^{- 1} \right)} \right)_{g, h \in G}}
$$
holds for any automorphism $\psi$ of a finite group $G$. 
\end{remark}

\begin{remark}\label{rem:4.4}
As mentioned in Section~$\ref{sec1.2}$, Ore~\cite{MR42365} obtained the explicit formula 
\begin{align*}
{\rm C}(x_{1}, x_{2}, \ldots, x_{n})  = \sum_{\lambda \in \Lambda_{n}^{1}} m_{\lambda}(\zeta_{(n, 1)}) x_{\lambda} = \sum_{\lambda \in \tilde{\Lambda}_{n}^{1}} m_{\lambda}(\zeta_{(n, 1)}) x_{\lambda}, 
\end{align*}
Brualdi and Newman~\cite{MR266850} proved 
$\Nu({\rm P}(\mathbb{Z} / n \mathbb{Z})) = |S(n, n)| = |\tilde{\Lambda}_{n}^{1}|$, 
and Thomas~\cite{MR2104820} showed that 
$\Nu(\Theta(\mathbb{Z} / n \mathbb{Z})) = \Nu({\rm P}(\mathbb{Z} / n \mathbb{Z}))$ holds for any prime power $n$. 
From these results and Remark~$\ref{rem:1.2}$, the following holds: 
Let $n$ be a prime power. Then, for any $\lambda \in \mathbb{Z}^{n}$, 
\begin{align}\label{eq:prime-power-nonvanishing}
|\lambda| \equiv 0 \!\pmod{n}
\quad \Longrightarrow \quad 
m_{\lambda}(\zeta_{(n, 1)}) \neq 0. 
\end{align}
Combining this with Remark~$\ref{rem:4.3}$~(2), we find that, for any prime power $n$ and $\lambda \in \mathbb{Z}^{n}$, 
$$ 
|\lambda| \not\equiv 0 \!\pmod{n}
\quad \Longleftrightarrow \quad 
m_{\lambda}(\zeta_{(n, 1)}) = 0. 
$$
\end{remark}

We give an elementary proof of \eqref{eq:prime-power-nonvanishing} in Remark~$\ref{rem:4.4}$ for the case when $n$ is a prime number. We use the following lemma. 

\begin{lem}\label{lem:4.5}
Let $\lambda = (\lambda_{1}, \lambda_{2}, \ldots, \lambda_{n}) \in \mathbb{Z}^{n}$ with $|\lambda| \equiv 0 \pmod{n}$ and let $f$ be a function of period $n$. 
Then 
$$
\sum_{\sigma \in \mathfrak{S}_{n}} f( \lambda_{1} \sigma(1) + \lambda_{2} \sigma(2) + \cdots + \lambda_{n} \sigma(n) ) 
= n \sum_{\tau \in \mathfrak{S}_{n - 1}} f( \lambda_{1} \tau(1) + \lambda_{2} \tau(2) + \cdots + \lambda_{n - 1} \tau(n - 1) ). 
$$
\end{lem} 
\begin{proof}
First, the following is true: 
$$
\sum_{\sigma \in \mathfrak{S}_{n}} f( \lambda_{1} \sigma(1) + \lambda_{2} \sigma(2) + \cdots + \lambda_{n} \sigma(n) )
= \sum_{i = 1}^{n} \sum_{\substack{\sigma \in \mathfrak{S}_{n} \\ \sigma(i) = n}} 
f( \lambda_{1} \sigma(1) + \lambda_{2} \sigma(2) + \cdots + \lambda_{n} \sigma(n) ). 
$$
Let $A_{i} := \sum_{\substack{\sigma \in \mathfrak{S}_{n} \\ \sigma(i) = n}} f( \lambda_{1} \sigma(1) + \lambda_{2} \sigma(2) + \cdots + \lambda_{n} \sigma(n) )$ for any $1 \leq i \leq n$. 
Proving 
$$
A_{1} = A_{2} = \cdots = A_{n} = \sum_{\tau \in \mathfrak{S}_{n - 1}} f(\lambda_{1} \tau(1) + \lambda_{2} \tau(2) + \cdots + \lambda_{n - 1} \tau(n - 1))
$$
is sufficient to complete the proof of the lemma. 
Since 
\begin{align*}
\left\{ \sigma(j) - \sigma(n) \mid 1 \leq j \leq n - 1 \right\} \equiv \left\{ 1, 2, \ldots, n - 1 \right\} \pmod{n}
\end{align*}
holds for any $\sigma \in \mathfrak{S}_{n}$, 
there uniquely exists $\tau_{\sigma} \in \mathfrak{S}_{n - 1}$ such that 
$$
\tau_{\sigma}(j) \equiv \sigma(j) - \sigma(n) \pmod{n}
$$
for any $1 \leq j \leq n - 1$. 
Therefore, the map $h_{i} \colon \{ \sigma \in \mathfrak{S}_{n} \mid \sigma(i) = n \} \ni \sigma \mapsto \tau_{\sigma} \in \mathfrak{S}_{n - 1}$ is well-defined. 
We prove $h_{i}$ is bijective. 
It is sufficient to show that $h_{i}$ is injective. 
If $h_{i}(\sigma) = h_{i}(\sigma')$, then 
$$
\sigma(j) - \sigma(n) \equiv \sigma'(j) - \sigma'(n) \pmod{n}
$$
for any $1 \leq j \leq n - 1$. 
\begin{enumerate}
\item[(i)] When $i = n$, from $\sigma(n) = \sigma'(n) = n$, we have $\sigma(j) \equiv \sigma'(j) \pmod{n}$ for any $1 \leq j \leq n$. 
\item[(ii)] When $i \neq n$, from $\sigma(i) = \sigma'(i) = n$, we have 
$$
- \sigma(n) \equiv \sigma(i) - \sigma(n) \equiv \sigma'(i) - \sigma'(n) \equiv - \sigma'(n) \pmod{n}. 
$$
This leads to $\sigma(j) \equiv \sigma'(j) \pmod{n}$ for any $1 \leq j \leq n$. 
\end{enumerate}
Therefore, $\sigma = \sigma'$. 
Thus, $h_{i}$ is bijective. 
Since 
\begin{align*}
\lambda_{1} \sigma(1) + \cdots + \lambda_{n} \sigma(n) &\equiv \lambda_{1} \sigma(1) + \cdots + \lambda_{n - 1} \sigma(n - 1) - (\lambda_{1} + \cdots + \lambda_{n - 1}) \sigma(n) \\ 
&\equiv \lambda_{1} \left\{ \sigma(1) - \sigma(n) \right\} + \cdots + \lambda_{n - 1} \left\{ \sigma(n - 1) - \sigma(n) \right\} \\ 
&\equiv \lambda_{1} \tau_{\sigma}(1) + \lambda_{2} \tau_{\sigma}(2) + \cdots + \lambda_{n - 1} \tau_{\sigma}(n - 1) \quad \pmod{n}
\end{align*}
holds and $h_{i}$ is bijective, 
we have 
$$
\sum_{\substack{\sigma \in \mathfrak{S}_{n} \\ \sigma(i) = n}} f( \lambda_{1} \sigma(1) + \cdots + \lambda_{n} \sigma(n) ) = \sum_{\tau \in \mathfrak{S}_{n - 1}} f(\lambda_{1} \tau(1) + \cdots + \lambda_{n - 1} \tau(n - 1))
$$
for any $1 \leq i \leq n$. 
This completes the proof. 
\end{proof}

\begin{thm}[Special case of \eqref{eq:prime-power-nonvanishing} in Remark~$\ref{rem:4.4}$]\label{thm:4.6}
Let $p$ be a prime number. 
Then, for any $\lambda \in \mathbb{Z}^{p}$, it holds that 
$$ 
|\lambda| \equiv 0 \!\pmod{p}
\quad \Longrightarrow \quad 
m_{\lambda}(\zeta_{(p, 1)})
\neq 0. 
$$
\end{thm}
\begin{proof}
Let $\lambda = (\lambda_{1}, \lambda_{2}, \ldots, \lambda_{p}) \in \mathbb{Z}^{p}$ 
and suppose that $|\lambda| \equiv 0 \pmod{p}$. 
From the expression~\eqref{eq:MSP} of the MSP, we have 
\[
m_{\lambda}(\zeta_{(p, 1)})
=
\left( \prod_{a \in \mathbb{Z}}\frac{1}{\lambda[a]!} \right)
\sum_{\sigma \in \mathfrak{S}_p}
\zeta_p^{\lambda_{1} \sigma(1) + \lambda_{2} \sigma(2) + \cdots + \lambda_{p} \sigma(p)}.
\]
Since the factor
\[
\prod_{a \in \mathbb{Z}}\frac{1}{\lambda[a]!}
\]
is nonzero, it remains to prove that 
\[
\sum_{\sigma \in \mathfrak{S}_p}
\zeta_p^{\lambda_{1} \sigma(1) + \lambda_{2} \sigma(2) + \cdots + \lambda_{p} \sigma(p)}
\neq 0.
\]
When $p = 2$, we can prove 
$\sum_{\sigma \in \mathfrak{S}_{2}} \zeta_{2}^{\lambda_{1}\sigma(1) + \lambda_{2}\sigma(2)} \neq 0$ by direct calculation. 
We may therefore assume that $p$ is an odd prime. 
From Lemma~$\ref{lem:4.5}$, we have 
\begin{align*}
\sum_{\sigma \in \mathfrak{S}_{p}} \zeta_{p}^{\lambda_{1}\sigma(1) + \lambda_{2}\sigma(2) + \cdots + \lambda_{p}\sigma(p)} 
&= p \sum_{\sigma \in \mathfrak{S}_{p-1}} \zeta_{p}^{\lambda_{1}\sigma(1) + \lambda_{2}\sigma(2) + \cdots + \lambda_{p-1}\sigma(p-1)}. 
\end{align*}
We prove $\sum_{\sigma \in \mathfrak{S}_{p-1}} \zeta_{p}^{\lambda_{1}\sigma(1) + \lambda_{2}\sigma(2) + \cdots + \lambda_{p-1}\sigma(p-1)} \neq 0$ by contradiction. 
For \(1\leq i\leq p\), define
\[
c_i
:=
| \left\{
\sigma\in \mathfrak{S}_{p-1}
\mid
\lambda_{1} \sigma(1) + \lambda_{2} \sigma(2) + \cdots + \lambda_{p-1} \sigma(p-1)
\equiv i \! \! \pmod p
\right\} |.  
\]
Then each $c_i$ is a non-negative integer. Moreover, every
$\sigma \in \mathfrak{S}_{p-1}$ contributes to exactly one of the numbers $c_1, c_{2}, \ldots, c_p$. 
Hence
\[
\sum_{i = 1}^{p} c_i = (p-1)!.
\]
By the definition of \(c_i\), 
we have
\[
\sum_{\sigma\in \mathfrak{S}_{p-1}}
\zeta_p^{\lambda_1\sigma(1)+\lambda_2\sigma(2)+\cdots+\lambda_{p-1}\sigma(p-1)}
=
\sum_{i=1}^{p} c_i \zeta_p^i.
\]
Assume now that $\sum_{i = 1}^{p} c_{i} \zeta_{p}^{i} = 0$. 
Then, since $\sum_{i = 1}^{p} \zeta_{p}^{i} = 0$ and $\{ \zeta_{p}, \zeta_{p}^{2}, \ldots, \zeta_{p}^{p-1} \}$ is linearly independent over $\mathbb{Q}$, 
it follows that $c_{1} = c_{2} = \cdots = c_{p}$.  
This implies that 
\[
p \mid \sum_{i = 1}^{p} c_{i} = (p-1)!. 
\]
This is a contradiction. 
Thus, we have $\sum_{\sigma \in \mathfrak{S}_{p-1}} \zeta_{p}^{\lambda_{1}\sigma(1) + \lambda_{2}\sigma(2) + \cdots + \lambda_{p-1}\sigma(p-1)} \neq 0$. 
\end{proof}

\clearpage

\bibliographystyle{plain}
\bibliography{reference}

\medskip
\begin{flushleft}
Naoya Yamaguchi\\
Faculty of Education \\ 
University of Miyazaki \\
1-1 Gakuen Kibanadai-nishi\\ 
Miyazaki 889-2192 \\
JAPAN\\
n-yamaguchi@miyazaki-u.ac.jp
\end{flushleft}

\medskip
\begin{flushleft}
Yuka Yamaguchi\\
Faculty of Education \\ 
University of Miyazaki \\
1-1 Gakuen Kibanadai-nishi\\ 
Miyazaki 889-2192 \\
JAPAN\\
y-yamaguchi@miyazaki-u.ac.jp
\end{flushleft}

\medskip
\begin{flushleft}
Genki Shibukawa\\
Faculty of Engineering\\ 
Kitami Institute of Technology\\
165 Koen-cho Kitami\\ 
Hokkaido 090-8507\\
JAPAN\\
g-shibukawa@mail.kitami-it.ac.jp
\end{flushleft}

\end{document}